\documentclass[12pt,fleqn]{scrartcl} 

\usepackage{framed}
\usepackage{hyperref}
\usepackage{doi}

\usepackage{amssymb}
\usepackage{amsmath}
\usepackage[T1]{fontenc}
\usepackage[latin1]{inputenc}
\usepackage{bm}

\usepackage{color}
\newcommand{\D}{\frak{D}}

\newcommand{\Q}{\hat{Q}}
\newcommand{\Laplace}{\bm \Delta}

\newcommand{\C}{C \kern -0.1em \ell}

\newcommand{\g}{\gamma}

\newcommand{\sqq}{\sqrt{q}} 
\newcommand{\vectorvar}{\underline{x}} 
\newcommand{\uX}{\underline{x}}

\usepackage{soul, cancel, xcolor}

\definecolor{dred}{RGB}{191 0 64}

\usepackage{ulem}

\newtheorem{thm}{Theorem}

\newtheorem{example}{Example}%
\newtheorem{rem}{Remark}%
\newtheorem{lem}{Lemma}
\newtheorem{prop}{Proposition}

\newtheorem{defn}{Definition}%

\begin{document}
	
\title{The $q$-Dirac Operator on  Quantum Euclidean Space}
\date{Swanhild Bernstein\footnote{\textsl{swanhild.bernstein@math.tu-freiberg.de}, Institute of Applied Analysis, TU Bergakademie Freiberg, Pr\"uferstr. 9, D-09599 Freiberg, Germany}, Martha Lina Zimmermann\footnote{\textsl{martha-lina.zimmermann@math.tu-freiberg.de}, Institute of Applied Analysis, TU Bergakademie Freiberg, Pr\"uferstr. 9, D-09599 Freiberg, Germany}, Baruch Schneider\footnote{\textsl{baruch.schneider@osu.cz}, University of Ostrava, 70103 Ostrava, Czech Republic}}
\maketitle
\pagenumbering{arabic} 
\begin{abstract} In this paper, we provide the foundations of quantum Clifford analysis in $q$-commu\-tative variables with symmetric difference operators. We consider a $q$-Dirac operator on the quantum Euclidean space that factorizes the $U_q(\frak{o})$-invariant Laplacian $\Delta_q.$ Due to the non-commutativity of the multiplication, we need a special Clifford algebra $\C_{0,n}^q.$ We define $q$-monogenic functions as null solutions of the $q$-Dirac operator and $q$-spherical monogenic functions. We define an inner Fischer product and decompose the space of homogeneous polynomials of degree $k.$
\end{abstract}
\section{Introduction}\label{sec1}
What is quantum calculus? First of all, in quantum calculus, we do not take limits. Derivatives are differences. Smoothness of functions is not required, even though the assumption of continuity makes sense in order to exclude erratic functions. The label ``quantum'' is used in very different contexts. One can discretize, replace commutative variables with $q$-commutative ones, and replace commutative algebras with non-commutative ones. \\[0.5ex]
$Q$-calculus has a long tradition and goes back to L. Euler in the 1740s when he initiated the theory of partitions, also called additive analytic number theory \cite{Euler1748, Ernst2012}. E. Heine, who studied with C. F. Gauß, P. G. L. Dirichlet, and C. G. J. Jacobi, propounded a theory for a so-called $q$-hypergeometric series. One of Heine's pupils was C. J. Thomae, who, together with reverend F. H. Jackson, would develop the so-called $q$-integral, the inverse to the $q$-derivative or $q$-difference operator. The derivative was invented by I. Newton and G. W. Leibniz. L. Euler and E. Heine had used variants of the $q$-derivative, but a real $q$-derivative was invented first by F. H. Jackson in 1908 \cite{Jackson1908}.
The theory developed as an alternative calculus and function theory. Specific calculi tried to overcome its shortcomings by using that many results are valid when rewritten using the analogs of $q$-calculus instead of the standard relations and denotations. That led to $q$-functions, special orthogonal $q$-polynomials, and Umbral calculus. 
\\[0.5ex]
The modern usage and application of $q$-calculus are in physics. Due to \cite{Fiore-Madore-2000} already, H. S. Snyder \cite{Snyder1947} suggested that the micro-structure of space-time at the Planck level might be better described using a non-commutative geometry. In 1954, W. Pauli suggested that the gravitational field might be considered as a universal regularizer for all quantum-field divergencies. One reason for using the notion ``quantum''  is that the $q$-differential is also a quantum differential, i.e., it fulfills the properties of a differential in quantum geometry. Secondly, quantum groups arose in the work of L.D. Faddeev and the Leningrad school on the inverse scattering method to solve integrable models. Around 1985, V.G. Drinfeld and M. Jimbo discovered a class of Hopf algebras, which can be considered as one-parameter deformations of universal enveloping algebras of semisimple complex Lie algebras, the so-called Drinfeld-Jimbo algebras. The invention of quantum groups is one of the outstanding achievements of mathematical physics and mathematics in the late 20th century \cite{KSch2011}. 
\\[0.5ex]
The Dirac operator is a physical and mathematical object. The Dirac equation describes the factorization of the D'Alembertian and, in a simplified version, the Laplacian. In an attempt to generalize complex function theory to higher dimensions, two theories are developed: harmonic analysis, where the Laplacian is the essential operator, and as a refinement of harmonic analysis, the so-called Clifford analysis based on a first-order differential operator, called Dirac-operator, that factorizes the Laplacian. This paper aims to define a $q$-deformed version of the Dirac operator related to the $q$-deformed Laplacian defined in \cite{NUW1996, IoKl2001}.
\\[0.5ex]
There are also attempts at $q$-Calculus in Clifford analysis. F. Sommen made the first attempt and found some similarities to the Clifford phase space \cite{Sommen1997} but needed to define an adequate  $q$-Dirac operator in the quantum plane. Unfortunately, the paper was forgotten, and a second attempt was made in \cite{CoSo2010} with commuting variables. The Dirac operator is defined axiomatically, mimicking the properties of the one-dimensional $q$ derivative. They proved that such an operator exists, but there is no explicit formula for the Dirac operator. 
\\[0.5ex]
In recent years, the topic has become popular again. There are approaches with commuting variables, and the Dirac operator is built up from partial $q$-derivatives \cite{BeZiSchn2022, GonzalezCervantes2024a, GonzalezCervantes2024}. We still believe an appropriate Dirac operator must be related to a $q$-Laplacian appropriately defined in the quantum plane. However, no quantum group would relate to $\Delta_q$ as the rotation group $SO(n)$ relates to the classical Laplace operator. The Drinfeld-Jimbo quantum algebras $U_q(\frak{so}_n)$ are not useful, because the reductions $U_q(\frak{so}_{n-1}) \subset U_q(\frak{so}_n)$ and the embedding $U_q(\frak{so}_n) \subset U_q(\frak{sl}_n)$ are not allowed. For this reason, an analog of the Gel'fand-Tsetlin basis for finite-dimensional irreducible representations of $U_q(\frak{so}_n)$ does not exist. But, there is a nonstandard $q$-deformation of the enveloping algebra $U^{\prime}_q(\frak{so}_n)$ of the Lie algebra $\frak{so}_n$ that gives the  new $q$-deformed associative algebra $U^\prime(\frak{so}_n).$ This algebra plays the role of the rotation group $SO(n)$ \cite{IoKl2001}. A similar observation was already made in \cite{NUW1996}. Our starting point is, therefore, the $U_q(\frak{o})$-invariant Laplacian defined in these papers. The definition of the quantum groups may be found in \cite{KSch2011} and the papers mentioned in the paper. 
\\[0.5ex]
The paper is organized as follows. After some denotations and a recap of Clifford algebras and Clifford analysis in Section 2 we introduce $q$-calculus for $q$-commuting variables in Section 3. Because Clifford analysis is a refinement of harmonic analysis, in Section 4 we start our considerations with $q$-harmonic functions and some of their properties. The factorization of the Laplacian leads to Dirac operators and the factorization of the $Q$-radius to vector variables. That will lead to different $q$-Clifford algebras. In Section 5, we will consider Clifford-valued polynomials in a $q$-Clifford algebra and define $q$-monogenic functions. 

\section{Preliminaries}
\subsection{Classical Clifford algebras}
We consider the Clifford algebra $\mathcal{C}\ell_{0,n}$ over $\mathbb{R}^n$ with the identity element $e_{\emptyset}$ satisfying $e_{\emptyset}^2 = 1$ and generating elements $e_1,...,e_n$ fulfilling the multiplication rules
\begin{equation*}
	e_ie_j + e_je_i = -2\delta_{ij}
	\label{eq:multrules}
\end{equation*}
for $i,j = 1,...,n$. The elements $e_i$ generate a basis $\mathcal{B}$ containing $2^n$ elements\\
$\mathcal{B} = (e_{\emptyset},e_1,e_2,...,e_1e_2,...,e_1e_2e_3,...)$. For $M:= \{1,...,n\}$ and $A:= \{(h_1,...,h_r \in \mathcal{P}M : 1 \leq h_1 \leq ... \leq h_r \leq n \}$ an arbitrary element of the Clifford algebra $\lambda \in \mathcal{C}\ell_{0,n}$ is written as
\begin{equation*}
	\lambda = \sum_A \lambda_A e_A\text{, } \lambda_A \in \mathbb{R}
\end{equation*}
with $e_A = e_{h_1}...e_{h_r}$. The length or norm of a Clifford number is
\begin{equation*}
	|\lambda |_0 := 2^{n/2}\left(\sum_{A} |\lambda_A|^2\right)^{1/2} .
\end{equation*}
\subsection{Clifford analysis}
We can define Clifford-valued functions $f$ in the same manner as Clifford numbers, i.e.
\begin{equation*}
	f(\underline{x}) =\sum_A f_A(\underline{x})e_A = \sum_A f_A(x_1,...,x_m)e_A\text{ with } f_A \colon \mathbb{R}^m \rightarrow \mathbb{R}.
\end{equation*}
A vector in $\mathbb{R}^n$ can be identified as vector variable $\underline{x}$
\begin{equation*}
	\underline{x} = \sum_{i=1}^n x_i e_i.
	\label{eq:vectorvar}
\end{equation*}
The multiplication rules imply that $\underline{x}^2$ is scalar valued, i.e. $\underline{x}^2 = -\sum_{i=1}^n x_i^2 = - \vert \underline{x} \vert^2$. Further, we can define the Dirac operator $\D_{\underline{x}}$
\begin{equation*}
	\D_{\underline{x}} =  \sum_{i=1}^n e_i \partial_i.
	\label{eq:Dirac}
\end{equation*}
Here $(\D_{\underline{x}})^2 = - \Delta$ is scalar with $\Delta$ denoting the Laplace operator. We can define the Euler operator $\mathbb{E}$ and the Gamma operator $\Gamma$ with these operators. The Euler operator follows from the anticommutator relation 
$$\{\D_{\underline{x}},\underline{x}\} = \D_{\underline{x}}\underline{x} + \underline{x}\D_{\underline{x}} = -2\mathbb{E}+n$$ and takes the form
\begin{equation*}
	\mathbb{E} = \sum_{i=1}^n x_i\partial_i.
\end{equation*}
Similarly, the Gamma operator can be obtained from the commutator relation $$[\underline{x},\D_{\underline{x}}] = \underline{x}\D_{\underline{x}} - \D_{\underline{x}} \underline{x} = 2\Gamma - n$$ and takes the form
\begin{equation*}
	\Gamma = \sum_{i<j} e_ie_j(x_i\partial_j - x_j\partial_i).
\end{equation*}
The vector variable and Dirac operator generate the Lie algebra $\mathfrak{osp}(1\vert 2)$. This is due to the following relations. 
\begin{prop}[\cite{BieXu2010}]
	The operators $\D_{\underline{x}}$ and $\underline{x}$ generate a Lie superalgebra, isomorphic to $\mathfrak{osp}(1\vert 2)$, with the following relations
	\begin{align*}
		&\{ \underline{x},\underline{x} \} = -2\vert \underline{x} \vert^2 \quad &\{ \D_{\underline{x}},\D_{\underline{x}}\} = -2\Laplace\\
		&\{\underline{x},\D_{\underline{x}}\} = -2\left(\mathbb{E}+\frac{n}{2}\right) \quad &\left[\mathbb{E}+\frac{n}{2},\D_{\underline{x}}\right] = -\D_{\underline{x}}\\
		&[\vert \underline{x} \vert^2,\D_{\underline{x}}] = -2\underline{x}\quad &\left[\mathbb{E}+\frac{n}{2},\underline{x}\right] = \underline{x}\\
		&[\Laplace,\underline{x}] = 2\D_{\underline{x}} \quad &\left[\mathbb{E}+\frac{n}{2},\Laplace\right] = -2\Laplace\\
		&[\Laplace^2,\vert\underline{x}\vert^2] = 4\left(\mathbb{E}+\frac{n}{2}\right) \quad &\left[\mathbb{E}+\frac{n}{2},\vert\underline{x}\vert^2\right] = 2\vert\underline{x}\vert^2,
	\end{align*}
	where $\mathbb{E}$ is the Euler operator.
\end{prop}
The even subalgebra, generated by $\D_{\underline{x}}^2$, $\underline{x}^2$ and $\mathbb{E}+\frac{n}{2}$, is isomorphic to $\mathfrak{sl}_2$ \cite{CoSo2010}. Further considerations lead to the Howe dual pair $(\mathrm{SO}(n), \mathfrak{sl}_2)$ \cite{Brackx2010}. This can be refined to the Howe dual pair $(\mathrm{Spin}(n),\mathfrak{osp}(1\vert 2))$ with $\mathrm{Spin}(n)$ being the double cover of $\mathrm{SO}(n)$.
Finally, a Clifford-valued function $f$ is called left (right) monogenic if and only if
\begin{equation*}
	\D_{\underline{x}}f = 0\quad (f\D_{\underline{x}} = 0).
\end{equation*}
If we look at the special case of monogeneous polynomials, it becomes apparent that monogeneous polynomials of degree $k$ are also homogeneous polynomials of degree $k$. With this, they are also eigenfunctions of the Euler operator. Take a homogeneous polynomial $P_k$ of degree $k$, then $\mathbb{E}P_k = kP_k$.
More on Clifford analysis can be found in the following monographs \cite{BDS1982, DSS1992, GM1991}.

\section{$Q$-Calculus}
For a multi-dimensional realisation on quantum vector space let $\mathcal{A} = \mathbb{C}[x_1,...,x_n]$ be a ring generated by the elements $x_1,\ldots,x_n$ satisfying the $q$-commutation rules $x_ix_j = qx_jx_i \iff x_jx_i = q^{-1}x_ix_j$ for $i<j$, where $0<q<\infty,\ q\not=1.$ \\[0.5ex]
The symmetric $q$-number is defined by
\begin{equation*}
	[m]_q = \frac{q^m-q^{-m}}{q-q^{-1}}.
\end{equation*}
Then we obtain the relation $\partial^q(x^m) = [m]_qx^{m-1}$ for the $q$-difference operator.
\begin{defn}
We define the number $\{a\}_q$ for some $a\in\mathbb{R}$ as
\begin{equation*}
	\{a\}_q = \frac{a-a^{-1}}{q-q^{-1}}.
\end{equation*}
\end{defn}
 Due to the $q$-commutativity, there are two different $q$-difference operators. We define the partial $q$-difference operator from the right by
\begin{equation*}
	\partial_i^R = (x^R_i)^{-1}\frac{\g_i-\g_i^{-1}}{q-q^{-1}}.
\end{equation*}
Here $x^R_i$ is the multiplication from the right with $x_i$, while $x^L_i$ is the multiplication from the left with $x_i$. 
Due to the $q$-commuting variables, we can introduce a second kind of $q$-difference operator using the left multiplication operator
\begin{equation*}
	\partial_i^L = (x^L_i)^{-1}\frac{\g_i-\g_i^{-1}}{q-q^{-1}}.
\end{equation*}
The operators $\g_i^{\pm 1}$ act on a function $p(x_1,\ldots,x_n)$ as
\begin{equation*}
	\g_i^{\pm 1}p(x_1,\ldots,x_n) = p(x_1,\ldots,x_{i-1},q^{\pm 1}x_i,x_{i+1},\ldots,x_n).
\end{equation*}
Furthermore,
\begin{align*}
	\omega_i & = \gamma_1^{-1}\cdots \gamma_{i-1}^{-1} \gamma_{i+1}\cdots \gamma_n \iff \omega_i^{-1} = \gamma_1\cdots \gamma_{i-1} \gamma_{i+1}^{-1}\cdots \gamma_n^{-1}  \\
	\text{and} & \quad \gamma = \gamma_1\gamma_2\ldots \gamma_n. 
\end{align*}

In general, these operators do not commute. Instead, we obtain the following relations:
\begin{itemize}\itemsep1ex
	\item[1)] $x_i^Lx_j^L = q x_j^Lx_i^L$ and $x_i^Rx_j^R = q^{-1} x_j^Rx_i^R,$ for $i<j,$
	\item[2)] $\partial_i^L\partial_j^L = q\partial_j^L\partial_i^L$ and $\partial_i^R\partial_j^R = q^{-1} \partial_j^R\partial_i^R,$ for $i<j,$
	\item[3)] $\partial_i^Rx_j^R = q x_j^R\partial_i^R $ and $\partial_i^Lx_j^L = q^{-1} x_j^L\partial_i^L ,$ for $i< j,$
	\item[4)] $\partial_i^R x_j^L = x_j^L\partial_i^R$ and $\partial_i^Lx_j^R = x_j^R\partial_i^L,$ for $i\not= j,$ 
	\item[5)] $\gamma_i x_j^R = q^{\delta_{ij}} x_j^R \gamma_i$  and $\gamma_i x_j^L = q^{\delta_{ij}} x_j^L \gamma_i,$ 
	\item[6)] $\gamma_i \partial_j^R = q^{-\delta_{ij}} \partial_j^R \gamma_i$  and $\gamma_i \partial_j^L = q^{-\delta_{ij}} \partial_j^L \gamma_i,$ 
	\item[7)] $x_i^R = x_i^L \omega_i ^{-1}  \iff x_i^L = x_i^R \omega_i,$
	\item[8)] $\partial_i^R = \partial_i^L \omega_i \iff \partial_i^L = \partial_i^R \omega_i^{-1}.$
\end{itemize}

With these basic relations we can prove more relations. 

\begin{lem}[Weyl relations]
We have the following Weyl relations:
\begin{align*}
	x_i^L\partial^L_j - q^{-1} \partial_j^L x_i^L = \left\{ \begin{array}{cl} 0&  i<j, \\ -q^{-1}\gamma^{-1}_i & i=j. \end{array} \right. \quad 
	x_i^R\partial^R_j - q \partial_j^R x_i^R = \left\{ \begin{array}{cl} 0&  i<j, \\ q\gamma_i & i=j. \end{array} \right.
\end{align*}
\end{lem}
Proof: To give the reader an idea how to obtain these relations we will prove the first Weyl relation. All operators act on polynomials $P(x)\in \mathcal{A}.$ For $i\not=j$ we obtain
\begin{align*}
x_i^L\partial_j^L & = x_i^L\left( \frac{1}{x_j^L} \frac{\left( \gamma_j - \gamma_j^{-1}\right)}{q-q^{-1}}\right) = q^{-1}\frac{1}{x_j^L}\frac{\left( x_i^L\gamma_j - x_i^L\gamma_j^{-1}\right)}{q-q^{-1}} 
  = q^{-1}\frac{1}{x_j^L}\frac{ \left( \gamma_j - \gamma_j^{-1}\right)}{q-q^{-1}} x_i^L \\ 
  & = q^{-1} \partial_j^Lx_i^L 
\end{align*}
and for $i=j$
\begin{align*}
	x_i^L\partial_i^L & = x_i^L\left( \frac{1}{x_i^L} \frac{\left( \gamma_i - \gamma_i^{-1}\right)}{q-q^{-1}}\right) = \frac{1}{x_i^L}\frac{\left( x_i^L\gamma_i - x_i^L\gamma_i^{-1}\right)}{q-q^{-1}} 
	= \frac{1}{x_i^L}\frac{ \left( q^{-1}\gamma_i - q\gamma_i^{-1}\right)}{q-q^{-1}} x_i^L \\ 
	& = \frac{1}{x_i^L}\frac{ \left( q^{-1}\gamma_i -  q^{-1} \gamma_i^{-1} +  q^{-1} \gamma_i^{-1} - q\gamma_i^{-1}\right)}{q-q^{-1}} x_i^L = q^{-1} \partial_i^Lx_i^L - \frac{1}{x_i^L} \gamma_i^{-1} x_i^L \\
	& = q^{-1} \partial_i^Lx_i^L - q^{-1}\gamma_i^{-1}
\end{align*}
\hfill $\square$
\begin{lem}[Mixed Weyl relations]
Furthermore, we have the  mixed relations:
\begin{align*}
	x_i^L\partial^R_j - \partial_j^R x_i^L =  \left\{ \begin{array}{cl} 0&  i\not=j, \\ -\left[\frac{1}{2}\right]_q(q^{1/2}\gamma_i + q^{-1/2} \gamma_i^{-1})\omega_i & i=j. \end{array} \right. \\
    x_i^R\partial^L_j -  \partial_j^L x_i^R = \left\{ \begin{array}{cl} 0&  i\not=j, \\ -\left[\frac{1}{2}\right]_q(q^{1/2}\gamma_i + q^{-1/2} \gamma_i^{-1})\omega^{-1}_i & i=j. \end{array} \right.
\end{align*}
It is easily seen that all relations give the classical relations in case $q=1.$
\end{lem}
Proof: We prove only the first relation. In this proof the definition of $\omega_i$ will become clear. For $i\not=j$ we obtain
\begin{align*}
x_i^L\partial_j^R & = x_i^L \frac{1}{x_j^R}\left(\frac{\gamma_j - \gamma_j^{-1}}{q-q^{-1}}\right) =  \frac{1}{x_j^R}\left(\frac{\gamma_j - \gamma_j^{-1}}{q-q^{-1}}\right) x_i^L = \partial_j^R x_i^L
\end{align*}
and for $i=j$:
\begin{align*}
x_i^L\partial_i^R - \partial_i^R x_i^L & = x_i^L \frac{1}{x_i^R} \left(\frac{\gamma_i - \gamma_i^{-1}}{q-q^{-1}}\right) - \frac{1}{x_i^R} \left(\frac{\gamma_i - \gamma_i^{-1}}{q-q^{-1}}\right) x_i^L \\
& = \frac{1}{x_i^R} \left(\frac{q^{-1}\gamma_i - q\gamma_i^{-1}}{q-q^{-1}}\right) x_i^L -\frac{1}{x_i^R} \left(\frac{\gamma_i - \gamma_i^{-1}}{q-q^{-1}}\right) x_i^L \\
& =  \frac{1}{x_i^R}\left(\frac{(q^{-1}-1)\gamma_i - (q-1) \gamma_i^{-1}}{q-q^{-1}}\right) x_i^L \\
& =\left(\frac{q(q^{-1}-1)\gamma_i - q^{-1}(q-1) \gamma_i^{-1}}{q-q^{-1}}\right)  \frac{1}{x_i^R} x_i^R \omega_i, \\
 \text{where we use that} & \text{$\frac{1}{x_i^R} x_i^L = \frac{1}{x_i^R} x_i^R \omega_i,$ because we have to bring $\frac{1}{x_i}$ and $x_i$ together,} \\
& =\left(\frac{q^{1/2}(q^{-1/2}-q^{1/2})\gamma_i - q^{-1/2}(q^{1/2}-q^{-1/2}) \gamma_i^{-1}}{q-q^{-1}}\right) \omega_i \\
& = - \left[\frac{1}{2}\right]_q (q^{1/2}\gamma_i + q^{-1/2}\gamma_i^{-1})\omega_i
\end{align*}
\hfill $\square$

\begin{lem}[Product rules for $\partial_i^R$]
We have the two product rules:
\begin{align*}
	\partial_i^R(fg) & = \gamma_i^{-1}(f)\partial_i^R(g) + \partial_i^R(f)\gamma(g), \\
	\partial_i^R(fg) & = \gamma_i(f)\partial_i^R(g) + \partial_i^R(f)\gamma^{-1}(g).
\end{align*}
\end{lem}
Proof: 
\begin{align*}
\partial_i^R(fg) & = \frac{1}{x_i^R}\left(\frac{\gamma_i(fg)-\gamma_i^{-1}(fg)}{q-q^{-1}}\right)  \\ 
&=\frac{1}{x_i^R}\left(\frac{\gamma_i(f)-\gamma_i^{-1}(f)}{q-q^{-1}}\right)\gamma_i(g) + \frac{1}{x_i^R}\gamma_i^{-1}(f)\left(\frac{\gamma_i(g)-\gamma_i^{-1}(g)}{q-q^{-1}}\right) \\
& = q\left(\frac{\gamma_i(f)-\gamma_i^{-1}(f)}{q-q^{-1}}\right)\gamma_i\left(g\frac{1}{x_i}\right) + \gamma_i^{-1}(f)\partial_i^R(g) 
 = \gamma_i^{-1}(f)\partial_i^R(g) + \partial_i^R(f)\omega_i^{-1}\gamma_i(g)  \\
& = \gamma_i^{-1}(f)\partial_i^R(g) + \partial_i^R(f)\gamma(g).
\end{align*}
Similarly,
\begin{align*}
	\partial_i^R(fg) & = \frac{1}{x_i^R} \gamma_i(f)\left(\frac{\gamma_i(g)-\gamma_i^{-1}(g)}{q-q^{-1}}\right) +  \frac{1}{x_i^R}\left(\frac{\gamma_i(f)-\gamma_i^{-1}(f)}{q-q^{-1}}\right) \gamma_i^{-1}(g) \\
	& = \gamma_i(f)\partial_i^R(g) + \partial_i^R(f)\omega_i\gamma_i^{-1}(g) = \gamma_i(f)\partial_i^R(g) +
	\partial_i^R(f)\gamma^{-1}(g).
\end{align*}
\hfill $\square$\\
After these basic relations, we recall some results for $q$-harmonic functions.

\section{Harmonic Analysis}
\subsection{$q$-Harmonic Analysis}
The theory of harmonic functions is well known. Specifically, the connection to group theory is studied by the group of N. Ja. Vilenkin. Here, we only mention the books \cite{Vilenkin1968, ViKl1993}. Later on, N.~Ja. Vilenkin and his group will also study $q$-harmonic functions.
M.~Noumi, T.~Umeda, and M.~Wakayama \cite{NUW1996} introduced a $q$-analogue of the Laplace operator on the considered quantum vector space that is invariant on $U_q(\mathfrak{o}_n)$
\begin{equation*}
	\Laplace^R_q = q^{n-1}(\partial_1^R)^2 + q^{n-2}(\partial_2^R)^2 + \ldots + (\partial_n^R)^2 = \sum_{i=1}^n q^{n-i}(\partial_i^R)^2.
	\label{eq:qLapalce}
\end{equation*}
The analog to the square radius $r^2= \sum_{j=1} x_j^2,$ the so-called $Q$-radius, where 
$$ Q = x_1^2 + q^{-1}x_2^2 + \ldots + q^{-n+1}x_n^2 $$
and the multiplication operator 
$$\hat{Q}^L =(x_1^L)^2 + q^{-1}(x_2^L)^2 + \ldots + q^{-n+1}(x_n^L)^2 $$ 
turn out to be helpful in the definition of harmonic polynomials. We need some relations between $\Laplace^R_q$ and $\hat{Q}^L$ \cite{NUW1996}
\begin{align}
	\Laplace^R_q(\hat{Q}^L)^k - (\hat{Q}^L)^k\Laplace^R_q = (\hat{Q}^L)^{k-1}[2k]_q\{q^{2k+n-2}\g^2\}_q \label{eq:Comm}  \\
	\Laplace^R_q(Q^k) = Q^{k-1}[2k]_q[2k+n-2]_q. \nonumber
\end{align}
It turns out that the operators $\hat{Q}^L, \Laplace^R_q $ and $q^n\g^2$ result in an algebra isomorphic to  $U_q(\mathfrak{sl}_2)$ \cite{NUW1996, KSch2011}, if we choose
$$ E = \frac{1}{\sqrt{[2]_q}} \Delta^R_q,\quad F= \frac{1}{\sqrt{[2]_q}} \hat{Q}^L \quad \text{and} \quad K= q^n\g^2 .$$
The first relations are easily verified by
$$ \g^2\g^{-2} = \g^{-2}\g^2 = 1, \quad \g^2\hat{Q}^L\g^{-2} = q^2\hat{Q}^L \quad \text{and}\quad \g^2\Laplace^R_q\g^{-2} = q^{-2}\Laplace^R_q . $$
The second relation can be obtained from (\ref{eq:Comm}) for $k=1$ or directly verified:
\begin{equation*}\label{eq:Comm2}
	\frac{1}{[2]_q}(\Laplace^R_q\hat{Q}^L - \hat{Q}°L\Laplace^R_q) = \{q^{n}\g^2\}_q = \frac{q^n\g^2 - q^{-n}\g^{-2}}{q-q^{-1}} . 
\end{equation*}

Unfortunately, there exists no quantum group that relates to $\Laplace^R_q$ as the rotation group $\mathrm{SO}(n)$ relates to the classical Laplace operator $\Laplace$ \cite{IoKl2001}. Instead, a non-standard (non-Drinfeld-Jimbo) $q$-deformation of the universal enveloping algebra $U'_q(\mathfrak{so}_n)$ of the Lie algebra $\mathfrak{so}_n$ will be used. This deformation is obtained by deforming the defining relations (Serr\'e relations) for the generating elements of $U'_q(\mathfrak{so}_n)$ \cite{IoKl2001}.\\

Based on the setting we have introduced so far, A.U. Klimyk and N. Z. Iorgov \cite{IoKl2001} define $q$-harmonic polynomials.
\begin{defn}
	A polynomial $p\in \mathcal{A}$ is a $q$-harmonic polynomial if $\Laplace^R_q p=0$. The space of all $q$-harmonic polynomials is denoted by $\mathcal{H}$, the space of $q$-harmonic homogeneous polynomials of degree $m$ by $\mathcal{H}_m.$ 
\end{defn}
Then the space of $q$-harmonic polynomials can be decomposed into $q$-harmonic homogeneous polynomials:
$$ \mathcal{H} = \bigoplus_{m=0}^{\infty} \mathcal{H}_m. $$
The space of homogeneous polynomials $\mathcal{P}_m$ of degree $m$ can be decomposed into the direct sum $\mathcal{P}_m = \mathcal{H}_m \oplus Q\mathcal{P}_{m-2}$ \cite{IoKl2001}. Repeatedly applying this decomposition leads to
\begin{equation*}
	\mathcal{P}_m = \bigoplus_{0\leq 2j\leq m} Q^j\mathcal{H}_{m-2j}.
\end{equation*}
Notice the similarity to the Fischer decomposition of harmonic polynomials in classical Clifford analysis. This implies that $Q$ is the appropriate generalisation of the squared radius $r^2$ used in the standard case.

\begin{example}
For an explicit formula for the harmonic polynomials, we will start with the case $n=1$. Then the $q$-harmonic polynomials are just of the form $ax+b$. 
\end{example}
\begin{example}[\cite{NUW1996}] \label{Ex2}
In the case, $n=2$, as in the classical case, the $q$-harmonics are two-dimensional (with an exception if $m=0$) and break into two irreducibles after extending the base field $\mathbb{R}$ by adjoining $\mathrm{i}$. Explicitly, we then get
\begin{equation*}
	\mathcal{H}_m = \mathbb{R}z_q^{m} \oplus \mathbb{R}\bar{z}_q^{m},
\end{equation*}
with $z_q^{m} = z_0z_1\ldots z_{m-1}$, $\bar{z}_q^{m} = \bar{z}_0\bar{z}_1\ldots\bar{z}_{m-1}$, $z_r = x_1+\mathrm{i}q^rx_2$ and $\bar{z}_r = x_1 - \mathrm{i}q^rx_2$.
This means they are of the form
\begin{align*}
	h_m(x_1,x_2) = a&(x_1+\mathrm{i}x_2)(x_1+\mathrm{i}qx_2)(x_1+\mathrm{i}q^2x_2)\cdots(x_1+\mathrm{i}q^{m-1}x_2)\\
	&+b(x_1-\mathrm{i}x_2)(x_1-\mathrm{i}qx_2)(x_1-\mathrm{i}q^2x_2)\cdots(x_1-\mathrm{i}q^{m-1}x_2)
\end{align*}
with $a,\, b\in \mathbb{R}$. For $q=1$ these are just the harmonic functions for $n=2.$ \\[0.5ex]
Proof: For $n=2$ we have $\gamma = \gamma_1\gamma_2$ and we easily see that $\gamma_1^{-1}z_k = (q^{-1}x_1 + iq^kx_2) = q^{-1}z_{k+1}$ and $\gamma_2z_k = (x_1 + iq^{k+1}x_2) = z_{k+1}.$ Hence,
\begin{align*}
\partial_1^Rz_q^m & = \sum_{k=0}^{m-1} \gamma_1^{-1}(z_0\cdots z_{k-1})\cdot \partial_1^Rz_k \cdot \gamma(z_{k+1}\cdots z_{m-1})
\end{align*}
\begin{align*}
& = \sum_{k=0}^{m-1} q^{-k}(z_1\cdots z_{k})\cdot q^{m-k-1} \cdot (z_{k+1}\cdots z_{m-1}) 
 = z_1\cdots z_m \sum_{k=0}^{m-1} q^{m-2k-1} = [m]_qz_1\cdots z_{m-1} .
\end{align*}
Similarly,
\begin{align*}
	& \partial_2^Rz_q^m  = \sum_{k=0}^{m-1} \gamma_1^{-1}(z_0\cdots z_{k-1})\cdot \partial_2^Rz_k \cdot \gamma(z_{k+1}\cdots z_{m-1}) \\
	& = \sum_{k=0}^{m-1} (z_1\cdots z_{k})\cdot iq^{k}q^{-(m-k-1)} \cdot (z_{k+1}\cdots z_{m-1}) 
	 =i z_1\cdots z_m \sum_{k=0}^{m-1} q^{-m+2k+1} = i[m]_qz_1\cdots z_{m-1}.
\end{align*}
From $\Delta_q^R = (q\partial_1^R - i\partial_2^R)(\partial_1^R + i\partial_2^R)$ and $(\partial_1^R + i\partial_2^R)z_q^m = 0$ we obtain 
$\Delta_q^R z_q^m = 0.$ \hfill $\square$
\end{example}

\begin{example} 
From the previous example we can easily obtain Clifford-valued harmonic functions. Our first generalisation identifies the complex numbers with the Clifford algebra $C\!\ell_{0,1},$. Still, it does not correspond to the factorisation of the Laplacian, which we want to use in the $n$-dimensional case. We introduce another variable $x_0$. For that we have to identify $z = x + iy$ with $x_0e_{\emptyset} + x_1e_1 \in C\!\ell_{0,1},$ where $e_{\emptyset}= 1$ and $e_1^2 = -1. $ Hence,  $$\Delta_q^R = (q\partial_0^R - e_1\partial_1^R)(\partial_0^R + e_1\partial_1^R).$$
The operators $q\partial_0^R - e_1\partial_1^R$ and $\partial_0^R + e_1\partial_1^R$ may be considered generalized Cauchy-Riemann operators. Obviously, $z_q^m = (x_0 + e_1x_1)(x_0 + e_1qx_1)\cdots (x_0 + e_1q^{m-1})$ and $\overline{z}_q^m = (x_0 - e_1x_1)(x_0 - e_1qx_1)\cdots (x_0 - e_1q^{m-1})$ are harmonic. 
That does not fit the factorization of the Laplacian we will obtain later on, where we use the Dirac operator and not the Cauchy-Riemann operator. \\
But, we can replace $\mathrm{i}$ by $e_2$  and obtain 
$$ \Delta_q^R = (q\partial_1^R -e_2\partial_2^R)(\partial_1^R + e_2 \partial_2^R) = q(\partial_1^R)^2 + (\partial_2^R)^2.$$
Hence, with $z_k = x_1 + e_2 q^kx_2$ the functions $z_q^m$ and with $\overline{z}_k = x_1 - e_2 q^kx_2$ the functions $\overline{z}_q^m$ are $q$-harmonic. \\
Although we use the Laplace operator here, which corresponds to our general Laplace operator for $n\geq 2$, the factorisation here also differs from that in Dirac operators, as we will introduce later. Unfortunately, we have not (yet) obtained this factorisation in the case $n>2$. In our opinion, this is because this factorisation corresponds to a complex theory and would, therefore, have to be represented with a complex Clifford algebra. In this article, however, we only consider real Clifford algebras. However, the case $n=2$ can be bent so that complex numbers $x+iy$ can be identified with $x_1 + e_2 x_2$. 
\end{example}

 For $n\geq 3$ we have the following recursion formula:
\begin{prop}[\cite{IoKl2001}]
	The $q$-harmonic polynomials have the following properties:\\
	 If $h_m(\mathbf{x})\in\mathcal{H}_m$, then $\tilde{h}_{m-1}(\mathbf{x}) := \gamma_n^{-1}\partial_n^qh_m(\mathbf{x}) \in \mathcal{H}_{m-1}$ and
	\begin{equation*}
		\hat{h}_{m+1}(\mathbf{x}):= h_m(\mathbf{x})x_n - \frac{Q\gamma_n^{-1}\partial_n^qh_m(\mathbf{x})}{[n+2m-2]_q}\in \mathcal{H}_{m+1}.
	\end{equation*}	
\end{prop}

\begin{example}\label{Ex4} Obviously, $x_n \in \mathcal{H}_1.$ Then 
	$$x_n^2- \frac{Q}{[n]_q} = x_n^2 - \frac{1}{[n]_q}\left(x_1^2 + q^{-1}x_2^2 + \ldots + q^{-n+1}x_n^2\right)\in \mathcal{H}_2, $$
	which we will prove directly. 
\begin{align*}
 \Delta_q^R\left(x_n^2- \frac{Q}{[n]_q}\right) & = [2]_q - \frac{[2]_q}{[n]_q}\left(q^{n-1}+ q^{n-3}+ \ldots + q^{-n+1}\right) \\
  & = [2]_q\left(1- \frac{q-q^{-1}}{q^n-q^{-n}}\left(q^{n-1}+ q^{n-3}+ \ldots + q^{-n+1}\right)\right) =0.
\end{align*}
\end{example}
\subsection{Factorization of the Laplacian and the $Q$-Radius}

To obtain a $q$-Dirac operator, we want to factorize the $U_q(\frak{o})$-invariant  $q$-Laplacian \cite{NUW1996, IoKl2001}
$$  \Delta^R_{q} = q^{(n-1)}(\partial^R_1)^2 + q^{(n-2)}(\partial^R_2)^2 + \ldots +(\partial^R_n)^2  = \sum_{i=1}^n q^{n-i}(\partial_i^R)^2 .$$
In general, to achieve a factorisation $(\D^R_q)^2 = -\Laplace_q$ the $q$-Dirac operator has to take the form
\begin{equation}
	\D^R_q = \sqq^{n-1}e_1\partial_1^R + \sqq^{n-2}e_2\partial_2^R + \ldots + e_n\partial_n^R = \sum_{i=1}^n \sqq^{n-i}e_i\partial_i^R.
\end{equation}
Due to the $q$-commutation rule $\partial_i^R\partial_j^R = q^{-1}\partial_j^R\partial_i^R$ for $i<j$ we need $q$-anti-commutation rules for our set of generators $e_i$. The multiplcation rules in our $q$-commuting case changes to $e_ie_j = -qe_je_i$ for $i<j$ and $e_i^2 = -1$.

For the $q$-vector variable, a factorization of the quantum squared radius $Q$ is needed.
Therefore, for the $q$-analogue of the relation $\vectorvar^2 = -r^2$ we define the $q$-vector variable
\begin{equation}
	\vectorvar_q = x_1e^+_1 + \sqq^{-1}x_2e^+_2 + \ldots + \sqq^{-n+1}x_ne^+_n = \sum_{i=1}^n \sqq^{-i+1}x_ie^+_i.
\end{equation}
The second set of generators $e^+_i$ behaves oppositely to the previously used generators $e_i$. This second set of generators $e^+_1,\ldots, e^+_n$ is governed by the multiplication rules
\begin{equation}
	e^+_ie^+_j = -q^{-1}e^+_je^+_i \text{ for } i<j \text{ and } (e^+_i)^2 = -1.
\end{equation}
We summarize these results in the following definition:

\begin{defn}
	The $q$-Clifford algebra $\mathcal{C}\ell_{0,n}^{q}, q\in \mathbb{R}_+$ over $\mathbb{R}^n$ is defined as the algebra with the identity element $e_{\emptyset}$ satisfying $e_{\emptyset}^2 = 1$ and generating elements $e_1,...,e_n$ fulfilling the multiplication rules
	\begin{align*}
		e_j^2 & = -1, \\
		e_ie_j + qe_je_i & = 0, \text{ for  } i< j, \quad  i,j = 1,...,n.
		\label{eq:multrules}
	\end{align*}
\end{defn}

\begin{rem}
For $q=1,$ we get the standard Clifford algebra $\C _{0,n}.$
\end{rem}


To return to our setting, we have two Laplacians and two $Q$-radii, distinguished by how the multiplication takes place, from the right or the left. Hence, we have two Dirac operators and vector variables:
\begin{align*}
	\D_q^R = \sqq^{n-1}e_1\partial_1^R + \sqq^{n-2}e_2\partial_2^R + \ldots + e_n\partial_n^R = \sum_{i=1}^n \sqq^{n-i}e_i\partial_i^R, \quad \D_q^R\D_q^R = -\Laplace^R,\\
	\vectorvar_q^R = \sqq^{n-1}e_1 x_1^R + \sqq^{n-2}e_2 x_2^R + \ldots + e_n x_n^R = \sum_{i=1}^n \sqq^{n-i}e_i x_i^R, \quad \vectorvar_q^R\vectorvar_q^R = -\Q^R.
\end{align*}
and 
\begin{align*}
	\D_q^L = \partial^L_1e^+_1 + \sqq^{-1}\partial^L_2e^+_2 + \ldots + \sqq^{-n+1}\partial^L_ne^+_n = \sum_{i=1}^n \sqq^{-i+1}\partial^L_ie^+_i, \quad \D_q^L\D_q^L = -\Laplace^L, \\
	\vectorvar_q^L = x^L_1e^+_1 + \sqq^{-1}x^L_2e^+_2 + \ldots + \sqq^{-n+1}x^L_ne^+_n = \sum_{i=1}^n \sqq^{-i+1}x^L_ie^+_i, \quad \vectorvar_q^L\vectorvar_q^L = -\Q^L.
\end{align*}

Because the left Laplacian and the right $Q$-radius fulfill the relations, we expect group relations will also be satisfied in the Clifford case. Unfortunately, we could not get such a relationship, but we got the following relations. \\[0.5ex]

\begin{thm}
	The operators $$E_i= \frac{1}{\sqrt{\sqrt{q}+ \sqrt{q}^{-1}}}\partial_i^Le_i, \quad F_i=\frac{1}{\sqrt{\sqrt{q}+ \sqrt{q}^{-1}}} x_i^Le_i,\quad  K_i = \sqrt{q}^{-1}\gamma_i^{-1}, \qquad i=1,2, \ldots n,$$ 
	build the algebra $\frak{osp}(n\!\mid\!2n).$
\end{thm}
Proof: We compute that the following relations are fulfilled:
\begin{align*}
	K_i E_i K_i^{-1} & = \gamma^{-1}_i\frac{1}{\sqrt{\sqrt{q}+ \sqrt{q}^{-1}}} \partial_i^Le_i \gamma_i = q\frac{1}{\sqrt{\sqrt{q}+ \sqrt{q}^{-1}}}\partial^L_ie_i = qE_i , \\
	K_i F_i K_i^{-1} &  = \frac{1}{\sqrt{\sqrt{q}+ \sqrt{q}^{-1}}}\gamma^{-1}_i x_i^Le_i \gamma_i = q^{-1} \frac{1}{\sqrt{\sqrt{q}+ \sqrt{q}^{-1}}}x_i^L e_i = q^{-1} F_i, \\
(\sqrt{q}+ \sqrt{q}^{-1}) &	(E_iF_i + F_iE_i) =  (-\partial^L_ix^L_i - x^L_i\partial^L_i) \\ 
 & = -   \left(x_i^L\frac{1}{x_i^L}\left(\frac{\gamma_i - \gamma_i^{-1}}{q-q^{-1}}\right) + \frac{1}{x_i^L}\left(\frac{\gamma_i - \gamma_i^{-1}}{q-q^{-1}}\right) x_i^L \right) \\
	& = -  \left(\frac{\gamma_i - \gamma_i^{-1}}{q-q^{-1}} + \frac{q\gamma_i - q^{-1}\gamma_i^{-1}}{q-q^{-1}} \right) = -  \left(\frac{(1+q)\gamma_i - (1+q^{-1})\gamma_i^{-1}}{q-q^{-1}}\right)  \\
	& = -  (\sqrt{q} + \sqrt{q}^{-1}) \left(\frac{\sqrt{q}\gamma_i - \sqrt{q}^{-1}\gamma_i^{-1}}{q-q^{-1}}\right) \\
\iff & E_iF_i + F_iE_i	 =\frac{\sqrt{q}^{-1} \gamma_i^{-1}- \sqrt{q} \gamma_i }{q-q^{-1}}  = \{ \sqrt{q}^{-1}\gamma^{-1}_i\}  = \frac{K_i - K_i^{-1}}{q-q^{-1}}. 
\end{align*}
\hfill $\square$

\begin{rem} 
    \begin{itemize}
    \item[~]
	\item The same is true for 
	$$E_i= \frac{1}{\sqrt{\sqrt{q}+ \sqrt{q}^{-1}}}\partial_i^Re_i, \quad F_i= \frac{1}{\sqrt{\sqrt{q}+ \sqrt{q}^{-1}}} x_i^Re_i,\quad  K_i = \sqrt{q}^{-1}\gamma_i^{-1}, \qquad i=1,2, \ldots n.$$ 
	\item Without the Clifford generating vectors, we have an algebra $\frak{osp}(n\!\mid\!2n)$ too: $$ E_i= \frac{1}{\sqrt{\sqrt{q}+ \sqrt{q}^{-1}}}x_i^{L/R}, \quad F_i= \frac{1}{\sqrt{\sqrt{q}+ \sqrt{q}^{-1}}} \partial_i^{L/R},\quad  K_i = \sqrt{q}\gamma_i,  \qquad i=1,2, \ldots n.$$ 
	\end{itemize}
\end{rem}


\section{Twisted Fischer Decomposition}
Because the variables and the left partial derivatives are $q$-commutative, we have a natural ring homomorphism over $\mathbb{C}$ as $x_i \mapsto \partial_i^R.$ In fact, we have some kind of conjugation among the $q$-commutative variables and the derivatives. To obtain a Fischer inner product, we have to combine the Clifford conjugation and the conjugation caused by the $q$-commutativity. 

Up to now, we just looked at specific polynomials. To obtain a basis we consider monomials, i.e. polynomials $\uX^{\alpha} = x_1^{\alpha_1}\cdots x_m^{\alpha_m},$ where $\alpha = (\alpha_1,\ldots ,\alpha_m) \in \mathbb{N}^m.$

Then $\{x^{\alpha}: |\alpha| = k\}$ is the standard basis of $\mathcal{P}_k,$ the set of all homogeneous polynomials of degree $k.$ 

Firstly, we introduce a Fischer inner product on the space of $q$-commutative variables.

\begin{defn}[$q$-Fischer inner product]
	Let $p_k^1, p^2_k \in \mathcal{P}_k$ then
    a Fischer inner product on $\mathcal{P}_k$ is defined by 
	$$ \langle p^1_k, p^2_k\rangle_k = \overline{p^1_k}(\partial_1^R, \ldots , \partial_n^R)(p^2_k(x_1, \ldots , x_n))^*\mid_{\mathbf{x}=0} =  p^1_k(\partial_n^R, \ldots , \partial_1^R)(p^2_k(x_1, \ldots, x_n))^*\mid_{x=0},  $$
where $p^1_k(\partial_n^R, \ldots , \partial_1^R)$ means the revised order of application of the derivatives and $^*$ the complex conjuation of the polynomial. \\[1ex]
The extension to Clifford-valued homogeneous polynomials $P_k^1, P^2_k \in \mathcal{P}(k, \C^q_{0,n})$ is given by 
$$ \langle P^1_k, P^2_k\rangle_k = [\overline{P^1_k}(\partial_1^R, \ldots , \partial_n^R)(P^2_k(x_1, \ldots , x_n))^*]_0\mid_{\mathbf{x}=0} = \sum_{|\alpha| = k} [\alpha]! q^{\sum\limits_{i<j} \alpha_i\alpha_j} [\overline{a}_{\alpha}^1 a_{\alpha}^2]_0 ,  $$
where $[a]_0$ denote the scalar part of $a\in \C_{0,n}^q.$
\end{defn}

\begin{rem}
The Fischer inner product is an inner product on the space  $\mathcal{P}_k$ of complex-valued homogeneous polynomials of degree $k.$ The extension of the Fischer inner product to Clifford-valued homogeneous polynomials of degree $k$ is a Clifford-valued inner product \cite{BDS1982}.
\end{rem}

\begin{example} For the monomials $x^{\alpha}\in \mathcal{P}_k$ we obtain (we write $\partial_i^{\alpha_i}$ instead of $(\partial_i^R)^{\alpha_i}$):
\begin{align*}
\langle x^{\alpha}, x^{\alpha} \rangle_k & = \partial_n^{\alpha_n}\cdots \partial_1^{\alpha_1}(x_1^{\alpha_1}\cdots x_n^{\alpha_n}) = q^{\sum\limits_{j=2}^n \alpha_j\alpha_1} \partial_n^{\alpha_n}\cdots \partial_2^{\alpha_2}[(\partial_1^{\alpha_1}x_1^{\alpha_1})x_2^{\alpha_2}\cdots x_n^{\alpha_n}] \\
 & = q^{\sum\limits_{j=2}^n \alpha_j\alpha_1} [\alpha_1]! \partial_n^{\alpha_n}\cdots \partial_2^{\alpha_2}(x_2^{\alpha_2}\cdots x_n^{\alpha_n}) \\ & = q^{\sum\limits_{j=2}^n \alpha_j\alpha_1} q^{\sum\limits_{j=3}^n \alpha_j\alpha_2} [\alpha_1]! [\alpha_2]! \partial_n^{\alpha_n}\cdots \partial_3^{\alpha_3}(x_3^{\alpha_3}\cdots x_n^{\alpha_n}) \\ 
 & = q^{\sum\limits_{i<j} \alpha_i\alpha_j} [\alpha]! 
\end{align*} 
Obviously, $\langle x^{\alpha}, x^{\beta} \rangle_k  = 0$ for $\alpha \not= \beta .$
\end{example}

With respect to the Fischer inner product, we consider the adjoint $a^{\dagger}$ for an operator $a$ defined by $\langle a^{\dagger}p^1, p^2 \rangle = \langle p^1, a p^2 \rangle.$ Obviously, $(ab)^{\dagger} = b^{\dagger}a^{\dagger}.$ Since the inner product is symmetric (or hermitian), the adjoint is involutory: $(a^{\dagger})^{\dagger} = a. $

\begin{example} Let $p_1\in \mathcal{P}_k$ and $p_2\in \mathcal{P}_{k+1},$
	then
\begin{align*}
\langle x_i^Rp_1(x), p_2(x) \rangle_{k+1} & = \langle p_1(x)x_i, p_2(x) \rangle_{k+1} = \partial_i^Rp_1(\partial^R) p_2^*(x)  
 = p_1(\partial^R) \partial_i^R\omega_i^{-1} p_2^*(x) \\ & = p_1(\partial^R) \partial_i^L p_2^*(x) = \langle p_1(x), \partial_i^L p_2(x) \rangle_k.
\end{align*}
This implies 
$$ \langle e_i x_i^Rp_1(x), p_2(x) \rangle_{k+1} =  \langle x_i^Rp_1(x), \overline{e}_ip_2(x) \rangle_{k+1} = \langle p_1(x), \partial_i^L \overline{e}_ip_2(x) \rangle_{k+1} $$
and hence 
$$({^L\vectorvar_q})^{\dagger} = (x_1^Le_1 + \sqq^{-1}x^L_2e_2 + \ldots + \sqq^{-n+1}x^L_ne_n)^{\dagger} = \sum_{i=1}^n (\sqq^{-i+1}x^L_ie_i)^{\dagger} =  \sqrt{q}^{(n-1)}\D^R, $$
where the vector variable ${^L\vectorvar_q}=x_1^Le_1 + \sqq^{-1}x^L_2e_2 + \ldots + \sqq^{-n+1}x^L_ne_n$ is build up with the multiplication from the left by the varialbes $x_i$ but lives in the algebra $\C_{0,n}^q,$ which in particular means that ${^L\vectorvar_q}$ does not factorize $\hat{Q}^L$! Furthermore, ${^L\vectorvar_q}\overline{{^L\vectorvar_q}}$ cannot be a scalar.
\end{example}



\begin{lem} We have the following formulas for the adjoints \vspace*{-2ex}\\
\begin{align*}
\gamma_i^{\dagger} & = \gamma_i, \quad (x_i^R)^{\dagger} = \partial_i^L,\quad (\partial_i^L)^{\dagger} = x_i^R \quad
	(x_i^L)^{\dagger} = \partial_i^R, \\[1ex]
	\text{and } & \text{hence }\quad (\Q^L)^{\dagger} = q^{-n+1}\Delta^R, \quad (\Delta^R)^{\dagger} = q^{n-1}\Q^L,\\[-6ex]
\end{align*}
\begin{align*}
	 (\frak{D}_q^R)^{\dagger} & = \sqrt{q}^{n-1}\overline{e}_1(\partial_1^R)^{\dagger}  + \ldots + \overline{e}_n(\partial_n^R)^{\dagger} = -\sqrt{q}^{n-1}x_1^L e_1- \ldots - x_n^L e_n \\ 
	  & =-\sqrt{q}^{-(n-1)}\ ^L\!\underline{x}_q,\ 
	 \text{where}\ ^L\underline{x}_q  = x_1^Le_1 + \sqrt{q}^{-1} x^L_2e_2 + \ldots + \sqrt{q}^{-n +1} x^L_ne_n,
	 \end{align*}
$$ \text{and}\ (^L\underline{x}_q)^{\dagger} = \sqrt{q}^{(n-1)} \frak{D}_q^R $$
\end{lem}

\begin{defn}[Monogenic polynomials] Let $R(\underline{x})$ be a Clifford-valued polynomial as an element of the algebra generated by $\{x_1, \ldots, x_n, e_1, \ldots , e_n\},$ where $e_1, \ldots, e_n$ are the generating elements of the Clifford algebra $\C_{0,n}^q, 0<q<\infty$. Then $R(\underline{x})$ is called right monogenic if it satisfies the Dirac equation $\D^R_qR_k(\underline{x})=0,$ where $\D^R_q= \sum_{j=1}^n \sqrt{q}^{n-i}\partial_j^R e_j. $ \\
With $\mathcal{M}_k$ we denote the space of all homogeneous right monogenic polynomials of degree $k.$
\end{defn}
\begin{example}
Because of $\D^R_q\D^R_q = - \Delta_q^R $ we obtain for $n=2$ from Example \ref{Ex2} that for $z_k = x_1 + e_1 q^kx_2$ and $z^m_q = z_0z_1\cdots z_{m-1}$ the Clifford-valued polynomials 
\begin{align*}
\D^R_q z^m_q = (\sqrt{q}\partial_1^R + e_1 \partial_2^R) z_q^m = (\sqrt{q}+e_1)[m]_q z_1\cdots z_{m-1}
\end{align*}
are right monogenic. 
\end{example}

\begin{defn}[Spherical monogenics] Let $P_k(x)$ be homogeneous of degree $k$ as an element of the algebra generated by $\{x_1, \ldots, x_n, e_1, \ldots , e_n\};$ then $P_k(x)$ is called spherical monogenic of degree $k$ if it satisfies the Dirac equation $\D^R_qP_k(x)=0.$ 
\end{defn}


\begin{thm}
	For all $R \in \mathcal{P}_k$ and $P \in \mathcal{P}_{k+1}$
	\begin{equation}
		\langle ^L\uX R,P\rangle_{k+1} = - \langle R,\, \D_q^R P\rangle_{k}.
	\end{equation}
	Hence, for $k\in\mathbb{N},$ the space of Clifford-valued homogeneous polynomials of degree $k,$ denoted by $\mathcal{P}_k$ can be decomposed as
	$$ \mathcal{P}_k = \mathcal{M}_k^R \oplus \ ^L\!\underline{x} \mathcal{P}_{k-1}.$$
	Further, the subspaces $\mathcal{M}_k^R$ and $\ ^L\!\underline{x} \mathcal{P}_{k-1}$ of $\mathcal{P}_k$ are orthogonal with respect to the Fischer inner product.
\end{thm}
Hence, we arrive at the Fischer decomposition
$$ \mathcal{P}_k = \sum_{s=0}^k \ (^L\!\underline{x})^s\mathcal{M}_{k-s}^R . $$

Proof:
	As $\mathcal{P}_k = {^L\uX}\mathcal{P}_{k-1} \oplus (^L\uX\mathcal{P}_{k-1})^{\bot}$ it suffices to proof $\mathcal{M}^R_k = (^L\uX\mathcal{P}_{k-1})^{\bot}$. For the first inclusion take any $R_{k-1} \in \mathcal{P}_{k-1}$ and $R_k \in \mathcal{P}_k$. Suppose now that
	\begin{equation*}
		\langle \underbrace{^L\uX R_{k-1}}_{\in ^L\uX\mathcal{P}_{k-1}},\underbrace{R_k}_{\in(^L\uX\mathcal{P}_{k-1})^{\bot}} \rangle_{k}= 0.
	\end{equation*}
	Then we have $\langle R_{k-1},\,\D^R R_k \rangle_{k_1} = 0$ for each $R_{k-1} \in \mathcal{P}_{k-1}$. We put $R_{k-1} = \,\D^R R_k$ and it follows that $\D^R R_k = 0$. Therefore $R_k \in \mathcal{M}_k^R$ and finally $(^L\uX\mathcal{P}_{k-1})^{\bot} \subset \mathcal{M}^R_k$.\\
	For the other inclusion take $P_k \in \mathcal{M}^R_k$. For each $R_{k-1} \in \mathcal{P}_{k-1}$
	\begin{align*}
		\langle ^L\uX R_{k-1},P_k\rangle_{k} &= -\langle R_{k-1},\underbrace{\D^R P_k}_{=0}\rangle_{k-1} = 0.
	\end{align*}
	As $^L\uX R_{k-1} \in ^L\uX\mathcal{P}_{k-1}$ we have $P_k \in (^L\uX\mathcal{P}_{k-1})^{\bot}$ and accordingly $\mathcal{M}^R_k \subset (^L\uX\mathcal{P}_{k-1})^{\bot}$. Therefore $\mathcal{M}^R_k = (^L\uX\mathcal{P}_{k-1})^{\bot}$. \hfill $\square$

The Fischer decomposition of the space of homogeneous polynomials $\mathcal{P}_k$ follows directly from 
\begin{equation*}
	\uX\mathcal{P}_{k-1} = \uX\mathcal{M}^q_{k-1} \oplus \uX^2\mathcal{P}_{k-2},\, \uX^2\mathcal{P}_{k-2} = \uX^2\mathcal{M}^q_{k-2} \oplus \uX^3\mathcal{P}_{k-3},\ldots\, .
\end{equation*}

\section{Summary and Outlook}
We worked with $q$-commutative variables on quantum Euclidean space. Because the Euclidean Laplacian is an $SO(n)$ invariant operator, looking for a similar Laplacian in quantum Euclidean space is natural.
Such an operator, the $U_q(\frak{o})$-invariant Laplacian $\Delta_q^R$, was used in \cite{NUW1996, IoKl2001}. Because the Drinfeld-Jimbo algebras are not useful for this consideration, new quantum enveloping algebras $U'_q(\frak{so}_n)$ are introduced in \cite{IoKl2001}. \\[0.5ex]
Associated to the Laplacian $\Delta_q^R$ is the $Q$-radius. Clifford analysis is a harmonic analysis refinement; the Dirac operator factorizes the Laplacian, and the vector variable $\underline{x}$ factorizes the square radius. It turns out that the Clifford algebra $\C_{0,n}$ does not allow any of these factorizations and it needs two different Clifford algebras $\C_{0,n}^q$ and $\C_{0,n}^{-q}$ to be able to factorize the Laplacian $-\Delta_q^R = \D_q^R\D_q^R $ and $-Q = \underline{x}_q^L\underline{x}_q^L.$ \\[0.5ex]
We consider the Dirac operator $\D_q^R$ in $\C_{0,n}^q$ and define a Fischer inner product based on the difference operators $\partial_i^R.$ It turns out that  $(\frak{D}_q^R)^{\dagger}  =-\sqrt{q}^{-(n-1)}\ ^L\!\underline{x}_q \not = \underline{x}^L_q,$ because it is impossible to factorize the $q$-Laplacian $\Delta_q^R$ and $Q^L$ in the same $q$-Clifford algebra. The Fischer inner product leads to the decomposition of homogeneous Clifford-valued polynomials of degree $k.$ \\[0.5ex]
The most interesting unsolved problem is whether it is possible to factorize the $q$-Laplacian and the $Q^L$ in the same Clifford algebra. We will study this problem further.

\end{document}